\documentclass[12pt]{article}

\usepackage{latexsym}
\usepackage{amsmath}
\usepackage{amssymb,amsfonts,amscd}
\usepackage{graphicx}

\newtheorem{theorem}{Theorem}[section]

\newtheorem{rem}[theorem]{Remark}

\newcommand{\bC}{{\mathbb C}}

\newcommand\qed{{\hspace*{\fill}$\Box$\vskip12pt plus 1pt}}

\newcommand\sF{{\mathcal F}}

\newcommand\sV{{\mathcal V}}
\newcommand\sW{{\mathcal W}}
\newcommand\sX{{\mathcal X}}

\newcommand\pn[1]{{\mathbb P}^{#1}}

\newcommand\W{W}
\newcommand\hatW{{\widehat\W}}
\newcommand\hatw{{\widehat\w}}

\newcommand\x{{\bf x}}
\newcommand\y{{\bf y}}
\newcommand\zero{{\bf 0}}

\newcommand\bfu{{\bf u}}
\newcommand\bfv{{\bf v}}
\newcommand\w{{\bf w}}

\newcommand{\Q}{\setminus Q}

\newcommand{\Witness}{{\bf Witness}}
\newcommand{\SysBySys}{{\bf SysBySys}}

\newcommand{\discard}{\includegraphics[bb = 0 0 326 309,width=0.7in]{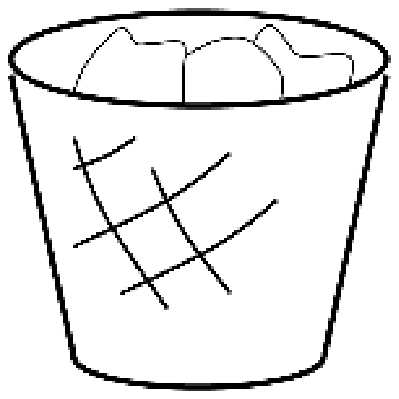}}
\newcommand{\YNboxLR}[4]{
  \begin{picture}(0,0)(0,0)
    \put(0,-3){\framebox(#1,6){#2}}
    \put(0,0){\makebox(0,0)[r]{#3$\,$}}
    \put(0,-3){\line(-5, 3){5}} \put(0,3){\line(-5,-3){5}}
    \put(#1,0){\makebox(0,0)[l]{$\,$#4}}
    \put(#1,-3){\line(5, 3){5}} \put(#1,3){\line(5,-3){5}}
  \end{picture}
}
\newcommand{\YNboxLB}[5]{
  \begin{picture}(0,0)(0,0)
    \put(0,-3){\framebox(#1,6){#2}}
    \put(0,0){\makebox(0,0)[r]{#3$\,$}}
    \put(0,-3){\line(-5, 3){5}} \put(0,3){\line(-5,-3){5}}
    \put(#4,-3.5){\makebox(0,0)[t]{$\,$#5}}
    \put(#4,-8){\line(3,5){3}} \put(#4,-8){\line(-3,5){3}}
  \end{picture}
}
\newcommand{\YNboxBR}[5]{
  \begin{picture}(0,0)(0,0)
    \put(0,-3){\framebox(#1,6){#2}}
    \put(#3,-3.5){\makebox(0,0)[t]{$\,$#4}}
    \put(#3,-8){\line(3,5){3}} \put(#3,-8){\line(-3,5){3}}
    \put(#1,0){\makebox(0,0)[l]{$\,$#5}}
    \put(#1,-3){\line(5, 3){5}} \put(#1,3){\line(5,-3){5}}
  \end{picture}
}
\newcommand{\WitnessArrayBox}[1]{
  \begin{picture}(0,0)(0,0)
    \put( 0,0){\framebox(30,6){}}
    \put( 5,3){\makebox(0,0){$\cdots$}}
    \put(15,3){\makebox(0,0){#1}}
    \put(25,3){\makebox(0,0){$\cdots$}}
    \put(10,0){\line(0,1){6}}
    \put(20,0){\line(0,1){6}}
  \end{picture}
}

\begin{document}

\title{Solving Polynomial Systems Equation by Equation\thanks{The
authors acknowledge the support of Land Baden-W\"urttemberg
(RiP-program at Oberwolfach).}}

\author{Andrew J. Sommese\thanks{Department of Mathematics,
University of Notre Dame, Notre Dame, IN 46556-4618, USA {\em
Email:} sommese@nd.edu {\em URL:} http://www.nd.edu/{\~{}}sommese.
This material is based upon work supported by the National Science
Foundation under Grant No.\ 0105653 and Grant No.\ 0410047; and
the Duncan Chair of the University of Notre Dame.}
\and Jan Verschelde\thanks{Department of Mathematics, Statistics,
and Computer Science, University of Illinois at Chicago, 851 South
Morgan (M/C 249), Chicago, IL 60607-7045, USA {\em Email:}
jan@math.uic.edu or jan.verschelde@na-net.ornl.gov {\em URL:}
http://www.math.uic.edu/{\~{}}jan. This material is based upon work
supported by the National Science Foundation under Grant No.\
0105739, Grant No.\ 0134611, and Grant No.\ 0410036.}
\and Charles W.  Wampler\thanks{General Motors Research and Development,
Mail Code 480-106-359, 30500 Mound Road, Warren, MI 48090-9055, U.S.A.
{\em Email:} Charles.W.Wampler@gm.com. This material is based upon work
supported by the National Science Foundation under Grant No.\
0410047 and by General Motors Research and Development.}}

\date{March 29, 2005}

\maketitle

\begin{abstract}
\noindent By a numerical continuation method called a diagonal homotopy
we can compute the intersection of two positive dimensional solution
sets of polynomial systems.  This paper proposes to use this diagonal
homotopy as the key step in a procedure to intersect general solution sets.
Of particular interest is the special case where one of the sets
is defined by a single polynomial equation.
This leads to an algorithm for finding a numerical representation
of the solution set of a system of polynomial equations introducing
the equations one-by-one.
Preliminary computational experiments show this approach can exploit
the special structure of a polynomial system, which
improves the performance of the path following algorithms.

\noindent {\bf 2000 Mathematics Subject Classification.} Primary
65H10; Secondary 13P05, 14Q99, 68W30.

\noindent {\bf Key words and phrases.}  Algebraic set,
component of solutions, diagonal homotopy, embedding,
equation-by-equation solver,
generic point, homotopy continuation, irreducible component,
numerical irreducible decomposition,
numerical algebraic geometry, path following,
polynomial system, witness point, witness set.

\end{abstract}

\newpage
\section{Introduction} \label{Sec:Introduction}

Homotopy continuation methods provide reliable and efficient numerical
algorithms to compute accurate approximations to all isolated solutions
of polynomial systems, see e.g.~\cite{Li03} for a recent survey.
As proposed in~\cite{SW}, we can approximate a positive dimensional
solution set of a polynomial system by isolated solutions, which
are obtained as intersection points of the set with a generic
linear space of complementary dimension.

New homotopy algorithms have been developed in a series of
papers~\cite{SV,SVW1,SVW4,SVW9,SVW10} to give numerical
representations of positive dimensional solution sets of polynomial
systems.  These homotopies are the main numerical algorithms in a
young field we call {\em numerical algebraic geometry}.
See~\cite{SW2} for a detailed treatment of this subject.

This paper provides an algorithm to compute numerical approximations
to positive dimensional solution sets of polynomial systems
by introducing the equations one at a time.
The advantage of working in this manner is that the special
properties of individual equations are revealed early in the
process, thus reducing the computational cost of later stages.
Consequently, although the new algorithm has more stages of
computation than earlier approaches, the amount of work in each stage
can be considerably less, producing a net savings in computing time.

This paper is organized in three parts.  First we explain our method
to represent and to compute a numerical irreducible decomposition of
the solution set of a polynomial system.  In the third section,
new diagonal homotopy algorithms will be applied to solve systems
subsystem by subsystem or equation by equation.  Computational
experiments are given in the fourth section.

\section{A Numerical Irreducible Decomposition}

We start this section with a motivating illustrative example,
which shows the occurrence of several solution sets, of different
dimensions and degrees.
Secondly, we define the notion of witness sets, which we
developed to represent pure dimensional solution sets of polynomial
systems {\em numerically}.  Witness sets are computed by
cascades of homotopies between embeddings of polynomial systems.

\subsection{An Illustrative Example}

Our running example (used also in~\cite{SVW1}) is
the following:

\begin{equation} \label{Eq:Illusex}
  f(x,y,z)
   = \left[
     \begin{array}{r}
        (y-x^2)(x^2+y^2+z^2-1)(x-0.5) \\
        (z-x^3)(x^2+y^2+z^2-1)(y-0.5) \\
        (y-x^2)(z-x^3)(x^2+y^2+z^2-1)(z-0.5) \\
     \end{array}
  \right].
\end{equation}
In this factored form we can easily identify the decomposition of the
solution set $Z=f^{-1}(\zero)$ into irreducible solution components,
as follows:
\begin{equation} \label{Eq:Illussol}
   Z = Z_2 \cup Z_1 \cup Z_0
     = \{Z_{21}\} \cup \{Z_{11} \cup Z_{12}  \cup Z_{13} \cup Z_{14} \}
       \cup \{Z_{01}\}
\end{equation}
where
\begin{center}
\begin{tabular}{l}
 1. $Z_{21}$ is the sphere $x^2+y^2+z^2-1=0$, \\
 2. $Z_{11}$ is the line $(x=0.5,z=0.5^3)$, \\
 3. $Z_{12}$ is the line $(x=\sqrt{0.5},y=0.5)$, \\
 4. $Z_{13}$ is the line $(x=-\sqrt{0.5},y=0.5)$, \\
 5. $Z_{14}$ is the twisted cubic $(y-x^2=0,z-x^3=0)$, \\
 6. $Z_{01}$ is the point $(x=0.5,y=0.5,z=0.5)$.
\end{tabular}
\end{center}

The sequence of homotopies in~\cite{SV} required to track 197 paths
to find a numerical representation of the solution set~$Z$.
With the new approach we will just have to trace 13 paths!
We show how this is done in Figure~\ref{Fig:Flowillusex}
in~\S \ref{Sec:Illusex} below, but we first describe a numerical
representation of~$Z$ in the next section.

\subsection{Witness Sets}

We define witness sets as follows. Let $f:\bC^N\rightarrow\bC^n$
define a system $f(\x) = \zero$ of $n$ polynomial equations $f =
\{f_1,f_2,\ldots,f_n \}$ in $N$ unknowns~$\x =
(x_1,x_2,\ldots,x_N)$. We denote the solution set of $f$ by
\begin{equation}
  V(f) = \{ \ \x \in\bC^N \ | \ f(\x)=\zero \ \}.
\end{equation}
This is a reduced\footnote{``Reduced'' means the set occurs with
multiplicity one, we ignore multiplicities $> 1$ in this paper.}
algebraic set. Suppose $X\subset V(f)\subset \bC^N$ is a pure
dimensional\footnote{``Pure dimensional'' (or ``equidimensional'')
means all components of the set have the same dimension.} algebraic
set of dimension~$i$ and degree~$d$. Then, a witness set for $X$ is
a data structure consisting of the system $f$, a generic linear
space $L\subset\bC^N$ of codimension~$i$, and the set of $d$ points
$X\cap L$.

If $X$ is not pure dimensional, then a witness set for $X$ breaks up
into a list of witness sets, one for each dimension. In our work, we
generally ignore multiplicities, so when a polynomial system has a
nonreduced solution component, we compute a witness set for the
reduction of the component. Just as $X$ has a unique decomposition
into irreducible components, a witness set for $X$ has a
decomposition into the corresponding irreducible witness sets,
represented by a partition of the witness set representation
for~$X$. We call this a \emph{numerical irreducible decomposition}
of~$X$.

The irreducible decomposition of the solution set~$Z$
in~(\ref{Eq:Illussol}) is represented by
\begin{equation}
  [W_2 , W_1, W_0] =
  [ [ W_{21} ], [ W_{11}, W_{12}, W_{13}, W_{14} ], [ W_{01} ] ],
\end{equation}
where the $W_i$ are witness sets for pure dimensional components,
of dimension~$i$, partitioned into witness sets $W_{ij}$'s corresponding
to the irreducible components of~$Z$.
In particular:
\begin{center}
\begin{tabular}{l}
 1. $W_{21}$ contains two points on the sphere,
    cut out by a random line, \\
 2. $W_{11}$ contains one point on the line $(x=0.5,z=0.5^3)$,
    cut out by a random plane, \\
 3. $W_{12}$ contains one point on the line $(x=\sqrt{0.5},y=0.5)$,
    cut out by a random plane, \\
 4. $W_{13}$ contains one point on the line $(x=-\sqrt{0.5},y=0.5)$,
    cut out by a random plane, \\
 5. $W_{14}$ contains three points on the twisted cubic,
    cut out by a random plane, \\
 6. $W_{01}$ is still just the point $(x=0.5,y=0.5,z=0.5)$.
\end{tabular}
\end{center}
Applying the formal definition, the witness sets $W_{ij}$ consist of
witness points $\w = \{ i, f, L, \x \}$, for $\x \in Z_{ij} \cap L$,
where $L$ is a random linear subspace of codimension~$i$ (in this
case, of dimension $3-i$). Moreover, observe $\#W_{ij} =
\deg(Z_{ij}) = \#(Z_{ij} \cap L)$.

Witness sets are set-theoretically equivalent to {\em lifting
fibers} which occur in a {\em geometric resolution} of polynomial
system. This geometric resolution is a symbolic analogue to a
numerical irreducible decomposition.  We refer to
\cite{GH93,GH01,GLS01,Lec03} for details about this symbolic
approach to solving polynomial system geometrically.

\subsection{Embeddings and Cascades of Homotopies}

A witness superset $\hatW_k$ for the pure $k$-dimensional part $X_k$
of~$X$ is a set in $X\cap L$, which contains $W_k:=X_k\cap L$ for a
generic linear space $L$ of codimension~$k$. The set of ``junk
points'' in $\hatW_k$ is the set $\hatW_k \setminus W_k$, which lies
in $ \left(\cup_{j>k}X_j\right)\cap L$.

The computation of a numerical irreducible decomposition for $X$
runs in three stages:
\begin{enumerate}
 \item Computation of a \emph{witness superset} $\hatW$ consisting
    of witness supersets $\hatW_k$ for each dimension $k = 1,2,\ldots,N$.
 \item Removal of junk points from $\hatW$ to get a witness set $W$ for $X$.

 \item Decomposition of $W$ into its irreducible components.
\newline In this stage, every witness set for a pure dimensional solution set
    is partitioned into witness sets corresponding to the irreducible
    components of the solution set.

\end{enumerate}

Up to this point, we have used the dimension of a component as the
subscript for its witness set, but in the algorithms that follow, it
will be more convenient to use codimension. The original algorithm
for constructing witness supersets was given in \cite{SW}. A more
efficient cascade algorithm for this was given in \cite{SV} by means
of an embedding theorem.

In \cite{SVW9}, we showed how to carry out the generalization of
\cite{SV} to solve a system of polynomials on a pure
$N$-dimensional algebraic set $Z\subset \bC^m$.  In the same
paper, we used this capability to address the situation where we
have two polynomial systems $f$ and $g$ on $\bC^N$ and we wish to
describe the irreducible decompositions of $A\cap B$ where
$A\in\bC^N$ is an irreducible component of $V(f)$ and $B\in\bC^N$
is an irreducible component of $V(g)$. We call the resulting
algorithm a \emph{diagonal homotopy}, because it works by
decomposing the diagonal system $\bfu-\bfv=\zero$ on $Z=A\times B$, where
$(\bfu,\bfv)\in\bC^{2N}$. In~\cite{SVW10}, we rewrote the homotopies
``intrinsically,'' which means that the linear slicing subspaces
are not described explicitly by linear equations vanishing on
them, but rather by linear parameterizations.
(Note that intrinsic forms were first used in a substantial way
to deal with numerical homotopies of parameterized linear spaces
in~\cite{HSS98}, see also~\cite{HV00}.)
This has always been allowed, even in \cite{SW},
but \cite{SVW10} showed how to do so consistently through the
cascade down dimensions of the diagonal homotopy, thereby
increasing efficiency by using fewer variables.

The subsequent steps of removing junk and decomposing the witness
sets into irreducible pieces have been studied in
\cite{SVW1,SVW2,SVW3,SVW4}.  These methods presume the capability to
track witness points on a component as the linear slicing space is
varied continuously.  This is straightforward for reduced solution
components, but the case of nonreduced components, treated in
\cite{SVW5}, is more difficult.  An extended discussion of the basic
theory may be found in \cite{SW2}.

In this paper, we use multiple applications of the diagonal homotopy
to numerically compute the irreducible decomposition of $A\cap B$
for general algebraic sets $A$ and $B$, without the restriction that
they be irreducible.  At first blush, this may seem an incremental
advance, basically consisting of organizing the requisite
bookkeeping without introducing any significantly new theoretical
constructs.  However, this approach becomes particularly interesting
when it is applied ``equation by equation,'' that is, when we
compute the irreducible decomposition of $V(f)$ for a system
$f=\{f_1,f_2,\ldots,f_n\}$ by systematically computing $V(f_1)$, then
$A_1\cap V(f_2)$ for $A_1$ a component of $V(f_1)$, then $A_2\cap
V(f_3)$ for $A_2$ a component of $A_1\cap V(f_2)$, etc. In this way,
we incrementally build up the irreducible decomposition one equation
at a time, by intersecting the associated hypersurface with all the
solution components of the preceding equations.  The main impact is
that the elimination of junk points and degenerate solutions at
early stages in the computation streamlines the subsequent stages.
Even though we use only the total degree of the equations---not
multihomogeneous degrees or Newton polytopes---the approach is
surprisingly effective for finding isolated solutions.

\section{Application of Diagonal Homotopies}

In this section, we define our new algorithms by means of two
flowcharts, one for solving subsystem-by-subsystem, and one that
specializes the first one to solving equation-by-equation. We then
briefly outline simplifications that apply in the case that only the
nonsingular solutions are wanted. First, though, we summarize the
notation used in the definition of the algorithms.

\subsection{Symbols used in the Algorithms}

 A witness set $W$ for a pure $i$-dimensional component $X$ in $V(f)$
is of the form $W=\{i,f,L,\sX\}$, where $L$ is the linear subspace
that cuts out the $\deg X$ points $\sX=X \cap L$.   In the following
algorithm, when we speak of a \emph{witness point} $\w\in W$, it
means that $\w=\{i,f,L,\x\}$ for some $\x\in\sX$. For such a $\w$
and for $g$ a polynomial (system) on $\bC^N$, we use the shorthand
$g(\w)$ to mean $g(\x)$, for $\x \in \w$.

In analogy to $V(f)$, which acts on a polynomial system, we
introduce the operator $\sV(W)$, which means the solution component
represented by the witness set~$W$.  We also use the same symbol
operating on a single witness point $\w=\{i,f,L,\x\}$, in which case
$\sV(\w)$ means the irreducible component of $V(f)$ on which point
$\x$ lies.  This is consistent in that $\sV(W)$ is the union of
$\sV(\w)$ for all $\w\in W$.

Another notational convenience is the operator $\sW(A)$, which gives
a witness set for an algebraic set $A$.  This is not unique, as it
depends on the choice of the linear subspaces that slice out the
witness points.  However, any two witness sets $W_1,W_2\in\sW(A)$
are equivalent under a homotopy that smoothly moves from one set of
slicing subspaces to the other, avoiding a proper algebraic subset
of the associated Grassmannian spaces, where witness points diverge
or cross.  That is, we have $\sV(\sW(A))=A$ and $\sW(\sV(W))\equiv W$,
where the equivalence in the second expression is under homotopy
continuation between linear subspaces.

The output of our algorithm is a collection of witness sets $W_i$,
$i=1,2,\ldots,N$, where $W_i$ is a witness set for the pure
\emph{codimension}~$i$ component of $V(f)$.  (This breaks from our
usual convention of subscripting by dimension, but for this
algorithm, the codimension is more convenient.) Breaking $W_i$ into
irreducible pieces is a post-processing task, done by techniques
described in \cite{SVW1,SVW3,SVW4}, which will not be described
here.

The algorithm allows the specification of an algebraic set
$Q\in\bC^N$ that we wish to ignore.  That is, we drop from the
output any components that are contained in $Q$, yielding witness
sets for $V(f_1,f_2,\ldots,f_n)\in\bC^N\Q$. Set $Q$ can be specified
as a collection of polynomials defining it or as a witness point
set.

For convenience, we list again the operators used in our notation,
as follows:
\begin{description}
    \item[$V(f)$] The solution set of $f(x)=0$.
    \item[$\sW(A)$] A witness set for an algebraic set $A$,
    multiplicities ignored, as always.
    \item[$\sV(W)$] The solution component represented by witness
    set $W$.
    \item[$\sV(\w)$] The irreducible component of $V(f)$ on which
    witness point $\w\in\sW(V(f))$ lies.
\end{description}

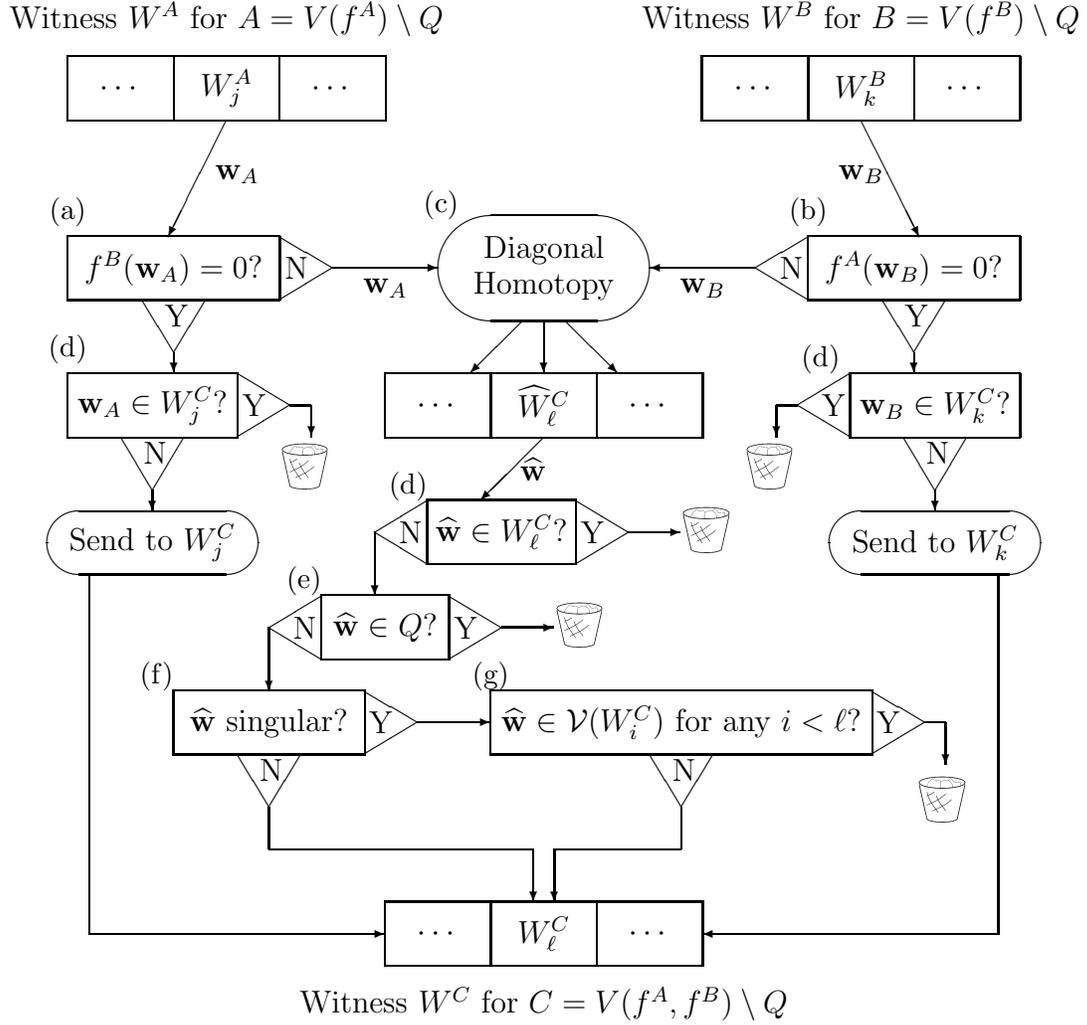
\begin{figure}
  \centering
  \setlength{\unitlength}{4pt}
  \def\sm{\small}
  \begin{picture}(92,110)(0,0)
     \put(15,100){\makebox(0,0)[b]{Witness $W^A$ for $A=V(f^A)\Q$}}
     \put(0,92){\WitnessArrayBox{$W^A_j$}}
     \put(75,100){\makebox(0,0)[b]{Witness $W^B$ for $B=V(f^B)\Q$}}
     \put(60,92){\WitnessArrayBox{$W^B_k$}}
     \put(30,12){\WitnessArrayBox{$W^C_\ell$}}
     \put(45,10){\makebox(0,0)[t]{Witness $W^C$ for $C=V(f^A,f^B)\Q$}}
     \put(15,92){\vector(-1,-2){5.5}}
     \put(16,86){\makebox(0,0)[b]{$\w_A$}}
     \put(0,78){\YNboxBR{20}{$f^B(\w_A)=0$?}{10}YN}
     \put(0,82){\makebox(0,0)[b]{\sm(a)}}
     \put(10,70){\vector(0,-1){2}}
     \put(75,92){\vector(1,-2){5.5}}
     \put(75,86){\makebox(0,0)[b]{$\w_B$}}
     \put(70,82){\makebox(0,0)[b]{\sm(b)}}
     \put(70,78){\YNboxLB{20}{$f^A(\w_B)=0$?}N{10}Y}
     \put(80,70){\vector(0,-1){2}}
     \put(45,78){\oval(20,10)} \put(37,84){\makebox(0,0)[r]{\sm(c)}}
     \put(45,78){\makebox(0,0){\shortstack{Diagonal\\ Homotopy}}}
     \put(25,78){\vector( 1,0){10}}
     \put(30,77){\makebox(0,0)[t]{$\w_A$}}
     \put(65,78){\vector(-1,0){10}}
     \put(60,77){\makebox(0,0)[t]{$\w_B$}}
     \put(45,73){\vector( 0,-1){5}} \put(43,73){\vector(-1,-1){5}}
     \put(47,73){\vector( 1,-1){5}}
     \put(30,62){\framebox(30,6){}}
     \put(35,65){\makebox(0,0){$\cdots$}}
     \put(45,65){\makebox(0,0){$\hatW^C_\ell$}}
     \put(55,65){\makebox(0,0){$\cdots$}} \put(40,62){\line(0,1){6}}
     \put(50,62){\line(0,1){6}}
     \put(0,65){\YNboxBR{16}{$\w_A\in W^C_j$?}8NY}
     \put(0,69){\makebox(0,0)[b]{\sm(d)}}
     \put(21,65){\line(1,0){2}}
     \put(23,65){\vector(0,-1){3}}
     \put(20,57){\discard}
     \put(8,57){\vector(0,-1){2}}
     \put(8,52){\oval(20,6)}
     \put(8,52){\makebox(0,0){Send to $W^C_j$}}
     \put(2,49){\line(0,-1){34}}
     \put(2,15){\vector(1,0){28}}
     \put(74,65){\YNboxLB{16}{$\w_B\in W^C_k$?}Y8N}
     \put(73,68){\makebox(0,0)[br]{\sm(d)}}
     \put(69,65){\line(-1,0){2}} \put(67,65){\vector(0,-1){3}}
     \put(64,57){\discard}
     \put(82,57){\vector(0,-1){2}}
     \put(82,52){\oval(20,6)}
     \put(82,52){\makebox(0,0){Send to $W^C_k$}}
     \put(88,49){\line(0,-1){34}}
     \put(88,15){\vector(-1,0){28}}
     \put(45,62){\vector(-1,-1){6}}
     \put(44,59){\makebox(0,0){$\hatw$}}
     \put(34,56){\makebox(0,0)[br]{\sm(d)}}
     \put(34,53){\YNboxLR{14}{$\hatw \in W^C_\ell$?}NY}
     \put(53,53){\vector(1,0){5}} \put(58,51){\discard}
     \put(29,53){\vector(0,-1){6}}
     \put(24,47){\makebox(0,0)[br]{\sm(e)}}
     \put(24,44){\YNboxLR{12}{$\hatw \in Q$?}NY}
     \put(41,44){\vector(1,0){5}}  \put(46,42){\discard}
     \put(19,44){\vector(0,-1){6}}
     \put(10,38){\makebox(0,0)[br]{\sm(f)}}
     \put(10,35){\YNboxBR{18}{$\hatw$ singular?}{9}NY}
     \put(19,27){\line(0,-1){4}} \put(19,23){\line(1,0){25}}
     \put(44,23){\vector(0,-1){5}}
     \put(33,35){\vector(1,0){7}} \put(40,38){\makebox(0,0)[b]{\sm(g)}}
     \put(40,35){\YNboxBR{36}{$\hatw \in \sV(W^C_i) \hbox{\ for any\ }i < \ell$?}{18}NY}
     \put(81,35){\line(1,0){2}} \put(83,35){\vector(0,-1){4}}
     \put(80.3,25.5){\discard}
     \put(58,27){\line(0,-1){4}}
     \put(58,23){\line(-1,0){12}} \put(46,23){\vector(0,-1){5}}
  \end{picture}
  \caption{Subsystem-by-subsystem generation of witness
  sets for \hbox{$V(f^A,f^B)\Q$.}}\label{Fig:SysBySys}
\end{figure}

\subsection{Solving Subsystem by Subsystem}

In this section, we describe how the diagonal homotopy can be
employed to generate a witness set $W=\sW(V(f^A,f^B)\Q)$, given
witness sets $W^A$ for $A=V(f^A)\Q$, and $W^B$ for $B=V(f^B)\Q$.
Let us denote this operation as $W=\SysBySys(A,B;Q)$.
Moreover, suppose
$\Witness(f;Q)$ computes a witness set $\sW(V(f)\Q)$ by any means
available, such as by working on the entire system $f$ as in our
previous works, \cite{SV,SW}, with junk points removed but not
necessarily decomposing the sets into irreducibles. With these two
operations in hand, one can approach the solution of any large
system of polynomials in stages.  For example, suppose
$f=\{f^A,f^B,f^C\}$ is a system of polynomials composed of three
subsystems, $f^A$, $f^B$, and $f^C$, each of which is a collection
of one or more polynomials. The computation of a witness set
$W=\sW(V(f)\Q)$ can be accomplished as
\begin{align*}
  W^A = \Witness(f^A;Q),&\quad W^B=\Witness(f^B;Q)\\
  W^C = \Witness(f^C;Q),&\quad W^{AB} = \SysBySys(W^A,W^B;Q),\\
  W = \SysBySys&(W^{AB},W^C;Q).
\end{align*}
This generalizes in an obvious way to any number of subsystems.
Although we could compute $W=\Witness(f;Q)$ by directly working on
the whole system $f$ in one stage, there can be advantages to
breaking the computation into smaller stages.

The diagonal homotopy as presented in \cite{SVW9} applies to
computing $A\cap B$ only when $A$ and $B$ are each irreducible.  To
implement \SysBySys, we need to handle sets that have more than one
irreducible piece. In simplest terms, the removal of the requirement
of irreducibility merely entails looping through all pairings of the
irreducible pieces of $A$ and $B$, followed by filtering to remove
from the output any set that is contained inside another set in the
output, or if two sets are equal, to eliminate the duplication.  In
addition to this, however, we would like to be able to proceed
without first decomposing $A$ and $B$ into irreducibles. With a bit
of attention to the details, this can be arranged.

Figure~\ref{Fig:SysBySys} gives a flowchart for algorithm \SysBySys.
For this to be valid as shown, we require that the linear subspaces
for slicing out witness sets are chosen once and for all and used in
all the runs of \Witness\ and \SysBySys.  In other words, the
slicing subspaces for $W^A$ and $W^B$ at the top of the algorithm
must be the same as each other and as the output~$W^C$. This ensures
that witness sets from one stage can, under certain circumstances,
pass directly through to the next stage. Otherwise, a continuation
step would need to be inserted to move from one slicing subspace to
another.

The setup of a diagonal homotopy to intersect two irreducibles
$A\in\bC^N$ and $B\in\bC^N$ involves the selection of certain random
elements.  We refer to \cite{SVW9,SVW10} for the full details.
All we need to know at present is that in choosing these random
elements the only dependence on $A$ and $B$ is their dimensions,
$\dim A$ and $\dim B$.  If we were to intersect another pair of
irreducibles, say $A'\in\bC^N$ and $B'\in\bC^N$, having the same
dimensions as the first pair, i.e., $\dim A'=\dim A$ and $\dim
B'=\dim B$, then we may use the same random elements for both.  In
fact, the random choices will be generic for any finite number of
intersection pairs.  Furthermore, if $A$ and $A'$ are irreducible
components of the solution set of the same system of polynomials,
$f^A$, and $B$ and $B'$ are similarly associated to system $f^B$,
then we may use exactly the same diagonal homotopy to compute $A\cap
B$ and $A'\cap B'$.  The only difference is that in the former case,
the start points of the homotopy are pairs of points
$(\alpha,\beta)\in \sW(A)\times\sW(B)\subset\bC^{2N}$, while in the
latter, the start points come from $\sW(A')\times\sW(B')$.

To explain this more explicitly, consider that the diagonal homotopy
for intersecting $A$ with $B$ works by decomposing $\bfu-\bfv$
on $A\times B$.
To set up the homotopy, we form the randomized system
\begin{equation}
  \sF(\bfu,\bfv) = \left[
  \begin{array}{l}
     R_A f^A(\bfu)\\
     R_B f^B(\bfv)
  \end{array}
  \right],
\end{equation}
where $R_A$ is a random matrix of size $(N-\dim A)\times \#(f^A)$
and $R_B$ is random of size $(N-\dim B)\times \#(f^B)$.  [By
$\#(f^A)$ we mean the number of polynomials in system $f^A$ and
similarly for $\#(f^B)$.]  The key property is that $A\times B$ is
an irreducible component of $V(\sF(\bfu,\bfv))$ for all $(R_A,R_B)$ in a
nonzero Zariski open subset of $\bC^{(N-\dim A)\times
\#(f^A)}\times\bC^{(N-\dim B)\times \#(f^B)}$, say $R_{AB}$. But
this property holds for $A'\times B'$ as well, on a possibly
different Zariski open subset, say $R_{A'B'}$.  But $R_{AB}\cap
R_{A'B'}$ is still a nonzero Zariski open subset, that is, almost
any choice of $(R_A,R_B)$ is satisfactory for computing both $A\cap
B$ and $A'\cap B'$, and by the same logic, for any finite number of
such intersecting pairs.

The upshot of this is that if we wish to intersect a
pure dimensional set $A=\{A_1,A_2\}\subset V(f^A)$ with a
pure dimensional set $B=\{B_1,B_2\}\subset V(f^B)$, where $A_1$,
$A_2$, $B_1$, and $B_2$ are all irreducible,  we may form one
diagonal homotopy to compute all four intersections $A_i\cap B_j$,
$i,j\in\{1,2\}$, feeding in start point pairs from all four
pairings.  In short, the algorithm is completely indifferent as to
whether $A$ and $B$ are irreducible or not.  Of course, it can
happen that the same irreducible component of $A\cap B$ can arise
from more than one pairing $A_i\cap B_j$, so we will need to take
steps to eliminate such duplications.

We are now ready to examine the details of the flowchart in
Figure~\ref{Fig:SysBySys} for computing $W^C=\sW(V(f^A,f^B)\Q)$ from
$W^A=\sW(V(f^A)\Q)$ and $W^B=\sW(V(f^B)\Q)$. It is assumed that the
linear slicing subspaces are the same for $W^A$, $W^B$, and $W^C$.
 The following items (a)--(g) refer to labels in that chart.
\begin{enumerate}
  \item[(a)] Witness point $\w_A$ is a generic point of the component of
  $V(f^A)$ on which it lies, $\sV(\w_A)$.
  Consequently, $f^B(\w_A)=0$ implies, with
  probability one, that $\sV(\w_A)$ is contained in some component
  of $V(f^B)$.  Moreover, we already know that $\w_A$ is not in any
  higher dimensional set of $A$, and therefore it cannot be in any
  higher dimensional set of $C$.  Accordingly, any point $\w_A$ that passes
  test~(a) is an isolated point in
  witness superset $\hatW^C$. The containment of $\sV(\w_A)$ in $B$ means
  that the dimension of the set is unchanged by intersection, so
  if $\w_A$ is drawn from $W^A_j$, its correct destination is $W^C_j$.

  On the other hand, if
  $f^B(\w_A)\ne0$, then $\w_A$ proceeds to the diagonal homotopy as
  part of the computation of $\sV(\w_A)\cap B$.
  \item[(b)] This is the symmetric operation to (a).
  \item[(c)] Witness points for components not completely contained in
  the opposing system are fed to the diagonal homotopy in order to
  find the intersection of those components.  For each combination
  $(a,b)$, where $a=\dim\sV(\w_A)$ and $b=\dim\sV(\w_B)$, there is a
  diagonal homotopy whose random constants are chosen once and for
  all at the start of the computation.
  \item[(d)] This test, which appears in three places, makes sure
  that multiple copies of a witness point do not make it into $W^C$.
  Such duplications can arise when $A$ and $B$ have components in
  common, when different pairs of irreducible components from $A$
  and $B$ share a common intersection component, or when some
  component is nonreduced.
  \item[(e)] Since a witness point $\hatw$ is sliced out generically from
  the irreducible component, $\sV(\hatw)$, on which it lies, if
  $\hatw \in Q$, then $\sV(\hatw)\subset Q$.  We have specified
  at the start that we wish to ignore such sets, so we throw them
  out here.
  \item[(f)] In this test, ``singular'' means that the Jacobian
  matrix of partial derivatives for the sliced system that cuts
  out the witness point is rank deficient.  We test this by a
  singular value decomposition of the matrix.  If the point is
  nonsingular, it must be isolated and so it is clearly a witness
  point.  On the other hand, if it is singular, it might be either
  a singular isolated point or it might be a junk point that lies
  on a higher dimensional solution set, so it must be subjected to
  further testing.
  \item[(g)] Our current test for whether a singular test point is
  isolated or not is to check it against all the higher dimensional
  sets.  If it is not in any of these, then it must be an isolated
  point, and we put it in the appropriate output bin.
\end{enumerate}

In the current state of the art, the test in box~(g) is done using
homotopy membership tests. This consists of following the paths of
the witness points of the higher dimensional set as its linear
slicing subspace is moved continuously to a generically disposed one
passing through the test point.  The test point is in the higher
dimensional set if, and only if, at the end of this continuation one
of these paths terminates at the test point, see~\cite{SVW2}.
In the future, it may
be possible that a reliable local test, based just on the local
behavior of the polynomial system, can be devised that determines if
a point is isolated or not.  This might substantially reduce the
computation required for the test.  As it stands, one must test the
point against all higher dimensional solution components, and so
points reaching box~(g) may have to wait there in limbo until all
higher dimensional components have been found.

The test~(e) for membership in $Q$ would entail a homotopy
membership test if $Q$ is given by a witness set.  If $Q$ is given
as $V(f^Q)$ for some polynomial system $f^Q$, then the test is
merely ``$f^Q(\hatw)=0?$'' We have cast the whole algorithm on
$\bC^N$, but it would be equivalent to cast it on complex projective
space $\pn N$ and use $Q$ as the hyperplane at infinity.

As a cautionary remark, note that the algorithm depends on $A$ and
$B$ being complete solution sets of the given polynomial subsystems,
excepting the same set $Q$.  It is not valid when $A$ or $B$ is a
partial list of components.  In particular, suppose $A$ and $B$ are
distinct irreducible components of the same system, i.e., $f^A=f^B$.
The diagonal homotopy applies to finding $A\cap B$, but if we feed
these into the current algorithm, we will not get the desired
result. This is because of tests (a) and~(b), which would pass the
witness points around the diagonal homotopy block and directly into
the output. The algorithm is designed to compute $V(f)$, which in
this case includes $A\cup B$.

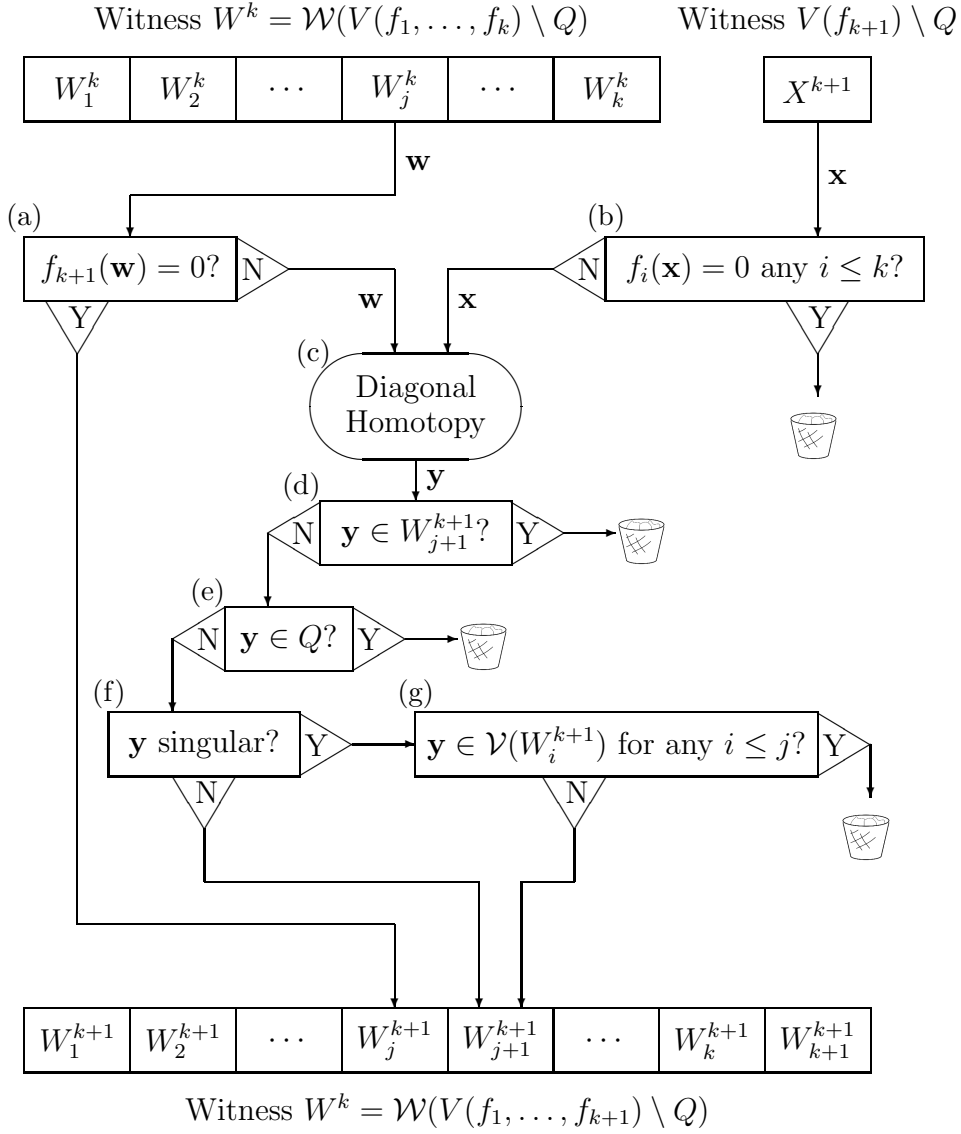
\begin{figure}
  \centering
  \setlength{\unitlength}{4pt}
  \def\sm{\small}
  \begin{picture}(90,120)(0,-10)
     \put(30,100){\makebox(0,0)[b]{Witness $W^k=\sW(V(f_1,\ldots,f_{k})\Q)$}}
     \put( 0,92){\framebox(60,6){}} \put( 5,95){\makebox(0,0){$W^k_1$}}
     \put(15,95){\makebox(0,0){$W^k_2$}}
     \put(25,95){\makebox(0,0){$\cdots$}}
     \put(35,95){\makebox(0,0){$W^k_j$}}
     \put(45,95){\makebox(0,0){$\cdots$}}
     \put(55,95){\makebox(0,0){$W^k_{k}$}} \put(10,92){\line(0,1){6}}
     \put(20,92){\line(0,1){6}} \put(30,92){\line(0,1){6}}
     \put(40,92){\line(0,1){6}} \put(50,92){\line(0,1){6}}
     \put(75,100){\makebox(0,0)[b]{Witness $V(f_{k+1})\Q$}}
     \put(70,92){\framebox(10,6){}}
     \put(75,95){\makebox(0,0){$X^{k+1}$}}
     \put( 0, 2){\framebox(80,6){}}
     \put( 5,5){\makebox(0,0){$W^{k+1}_1$}}
     \put(15,5){\makebox(0,0){$W^{k+1}_2$}}
     \put(25,5){\makebox(0,0){$\cdots$}}
     \put(35,5){\makebox(0,0){$W^{k+1}_j$}}
     \put(45,5){\makebox(0,0){$W^{k+1}_{j+1}$}}
     \put(55,5){\makebox(0,0){$\cdots$}}
     \put(65,5){\makebox(0,0){$W^{k+1}_{k}$}}
     \put(75,5){\makebox(0,0){$W^{k+1}_{k+1}$}}
     \put(10,2){\line(0,1){6}} \put(20,2){\line(0,1){6}}
     \put(30,2){\line(0,1){6}} \put(40,2){\line(0,1){6}}
     \put(50,2){\line(0,1){6}} \put(60,2){\line(0,1){6}}
     \put(70,2){\line(0,1){6}}
     \put(40,0){\makebox(0,0)[t]{Witness $W^k=\sW(V(f_1,\ldots,f_{k+1})\Q)$}}
     \put(35,92){\line(0,-1){7}} \put(35,85){\line(-1,0){25}}
     \put(10,85){\vector(0,-1){4}} \put(36,88){\makebox(0,0)[l]{$\w$}}
     \put(0,81.5){\makebox(0,0)[b]{\sm(a)}}
     \put(0,78){\YNboxBR{20}{$f_{k+1}(\w)=0$?}{5}YN}
     \put(5,70){\line(0,-1){54}} \put(5,16){\line(1,0){30}}
     \put(35,16){\vector(0,-1){8}}
     \put(75,92){\vector(0,-1){11}} \put(76,86){\makebox(0,0)[lb]{$\x$}}
     \put(55,81.5){\makebox(0,0)[b]{\sm(b)}}
     \put(55,78){\YNboxLB{30}{$f_i(\x)=0$ any $i\le k$?}N{20}Y}
     \put(75,70){\vector(0,-1){4}}
     \put(72.2 ,60){\discard}
     \put(37,65){\oval(20,10)} \put(29,70){\makebox(0,0)[r]{\sm(c)}}
     \put(37,65){\makebox(0,0){\shortstack{Diagonal\\ Homotopy}}}
     \put(25,78){\line(1,0){10}} \put(35,78){\vector(0,-1){8}}
     \put(50,78){\line(-1,0){10}} \put(40,78){\vector(0,-1){8}}
     \put(34,75){\makebox(0,0)[r]{$\w$}}
     \put(41,75){\makebox(0,0)[l]{$\x$}} \put(37,60){\vector(0,-1){4}}
     \put(38,58){\makebox(0,0)[l]{$\y$}}
     \put(28,56){\makebox(0,0)[br]{\sm(d)}}
     \put(28,53){\YNboxLR{18}{$\y\in W^{k+1}_{j+1}$?}NY}
     \put(51,53){\vector(1,0){5}}
     \put(56,50){\discard}
     \put(23,53){\vector(0,-1){7}}
     \put(19,46){\makebox(0,0)[br]{\sm(e)}}
     \put(19,43){\YNboxLR{12}{$\y\in Q$?}NY}
     \put(36,43){\vector(1,0){5}}
     \put(41,40){\discard}
     \put(14,43){\vector(0,-1){7}}
     \put( 8,36.5){\makebox(0,0)[b]{\sm(f)}}
     \put( 8,33){\YNboxBR{18}{$\y$ singular?}9NY}
     \put(31,33){\vector(1,0){6}}
     \put(17,25){\line(0,-1){5}}
     \put(17,20){\line(1,0){26}} \put(43,20){\vector(0,-1){12}}
     \put(37,36.5){\makebox(0,0)[b]{\sm(g)}}
     \put(37,33){\YNboxBR{38}{$\y \in \sV(W^{k+1}_i)
       \hbox{\ for any\ }i\le j$?}{15}NY}
     \put(80,33){\vector(0,-1){5}}
     \put(77.2,22){\discard}
     \put(52,25){\line(0,-1){5}}
     \put(52,20){\line(-1,0){5}}
     \put(47,20){\vector(0,-1){12}}
  \end{picture}
  \caption{Stage~$k$ of equation-by-equation generation of witness
  sets for $V(f_1,\ldots,f_n)\in\bC^N\Q$}\label{Fig:EqByEq}
\end{figure}

\subsection{Solving Equation by Equation}

The equation-by-equation approach to solving a polynomial system is
a limiting case of the subsystem-by-subsystem approach, wherein one
subsystem is just a single polynomial equation.  Accordingly, we
begin by computing a witness set $X^i$ for the solution set
$V(f_i)$, $i=1,2,\ldots,n$ of each individual polynomial.  If any
polynomial is identically zero, we drop it and decrement~$n$. If any
polynomial is constant, we terminate immediately, returning a null
result.  Otherwise, we find $X^i=(V(f_i)\cap L)\Q$, where $L$ is a
1-dimensional generic affine linear subspace.  A linear
parameterization of $L$ involves just one variable, so $X^i$ can be
found with any method for solving a polynomial in one variable,
discarding any points that fall in~$Q$.

Next, we randomly choose the affine linear subspaces that will cut
out the witness sets for any lower dimensional components that
appear in succeeding intersections.

The algorithm proceeds by setting $W^1=X^1$ and then computing
$W^{k+1}$ $=$ $\SysBySys$ $(W^k,X^{k+1};Q)$ for $k=1,2,\ldots,n-1$.  The
output of stage~$k$ is a collection of witness sets $W^{k+1}_i$ for
$i$ in the range from 1 to $\min(N,k+1)$.  (Recall, we are using the
codimension for the subscript.) Of course, some of these may be
empty, in fact, in the case of a total intersection, only the lowest
dimensional one, $W^{k+1}_{k+1}$, is nontrivial.

In applying the subsystem-by-subsystem method to this special case,
we can streamline the flowchart a bit, due to the fact that
$V(f_{k+1})$ is a hypersurface.  The difference comes in the
shortcuts that allow some witness points to avoid the diagonal
homotopy.

The first difference is at the output of test~(a), which now sends
$\w$ directly to the final output without any testing for duplicates.
This is valid because we assume that on input $\sV(\w)$ is not
contained within any higher dimensional component of $\sV(W^k)$, and
in the intersection with hypersurface $V(f_{k+1})$ that is the only
place a duplication could have come from.

On the opposing side, test~(b) is now stronger than before.  The
witness point $\x$ only has to satisfy one polynomial among
$f_1,f_2,\ldots,f_k$ in order to receive special treatment.  This is
because we already threw out any polynomials that are identically
zero, so if $f_j(\x)=0$ it implies that $\sV(\x)$ is a factor of
$V(f_j)$.  But the intersection of that factor with all the other
$V(f_i)$, $i\ne j$, is already in $V(f_1,f_2,\ldots,f_k)$, so nothing
new can come out of intersecting $\sV(\x)$ with $V(f_1,f_2,\ldots,f_k)$.
Accordingly, we may discard $\x$ immediately.

Another small difference from the more general algorithm is that the
test for junk at box~(g) never has to wait for higher dimensional
computations to complete.  When carrying out the algorithm, we draw
witness points from $W^k$ in order proceeding from left to right so
that computations are performed by decreasing dimension. Moreover,
we should run all the witness points in $W^k_j$ through test~(a)
before proceeding to feed any of them to the diagonal homotopy. This
ensures that all higher dimensional sets are in place before we
begin computations on $W^k_{j+1}$.  This is not a matter of much
importance, but it can simplify coding of the algorithm.

In the test at box~(d), we discard duplications of components,
including points that appear with multiplicity due to the presence
of nonreduced components. However, for the purpose of subsequently
breaking the witness set into irreducible components, it can be
useful to record the number of times each root appears. By the
abstract embedding theorem of~\cite{SVW9}, points on the same
irreducible component must appear the same number times, even though
we cannot conclude from this anything about the actual multiplicity
of the point as a solution of the system $\{f_1,f_2,\ldots,f_n\}$.
Having the points partially partitioned into subsets known to
represent distinct components will speed up the decomposition phase.

A final minor point of efficiency is that if $n>N$, we may arrive at
stage $k\ge N$ with some zero dimensional components, $W_N^k$. These
do not proceed to the diagonal homotopy: if such a point fails
test~(b), it is not a solution to system $\{f_1,f_2,\ldots,f_{k+1}\}=0$,
and it is discarded.

\subsection{Seeking only Nonsingular Solutions}

In the special case that $n\le N$, we may seek only the
multiplicity-one components of codimension $n$. (For $n=N$, this
means we seek only the nonsingular solutions of the system.) In this
case, we discard points that pass test~(a), since they give
higher dimensional components. Furthermore, we keep only the points
that test~(e) finds to be nonsingular and discard the singular ones.
This can greatly reduce the computation for some systems.

In this way, we may use the diagonal homotopy to compute nonsingular
roots equation-by-equation.  This performs differently than more
traditional approaches based on continuation, which solve the entire
system all at once.  In order to eliminate solution paths leading to
infinity, these traditional approaches use multihomogeneous
formulations or toric varieties to compactify $\bC^N$.  But this
does not capture other kinds of structure that give rise to positive
dimensional components.  The equation-by-equation approach has the
potential to expose some of these components early on, while the
number of intrinsic variables is still small, and achieves
efficiency by discarding them at an early stage. However, it does
have the disadvantage of proceeding in multiple stages.  For
example, in the case that all solutions are finite and nonsingular,
there is nothing to discard, and the equation-by-equation approach
will be less efficient than a one-shot approach. However, many
polynomial system of practical interest have special structures, so
the equation-by-equation approach may be commendable. It is too
early to tell yet, as our experience applying this new algorithm on
practical problems is very limited.  Experiences with some simple
examples are reported in the next section.

\section{Computational Experiments}\label{Sec:Examples}

The diagonal homotopies are implemented in the software
package PHCpack~\cite{V99}.  See~\cite{SVW7} for a description
of a recent upgrade of this package to deal with positive
dimensional solution components.

\subsection{An illustrative example}\label{Sec:Illusex}

\begin{figure}
\begin{center}
\begin{picture}(400,380)(-20,0)

   \put( 0,360){\framebox(80,20)[c]{$\#X^1 = 5$}}

   \put(33,353){\vector(0,-1){26}}  \put(35,343){${}_5$}

   \put(  0,300){\framebox(80,20)[c]{$\#W^1_1 = 5$}}
   \put( 80,300){\framebox(80,20)[c]{$W^1_2 = \emptyset$}}
   \put(160,300){\framebox(80,20)[c]{$W^1_3 = \emptyset$}}

   \put(33,293){\vector(0,-1){24}}  \put(35,283){${}_5$}
   \put(33,203){\vector(0,-1){26}}  \put(35,193){${}_2$}

   \put(0,206){
   \begin{picture}(60,60)(0,0)
   \put( 0,20){\framebox(60,40)[c]{$f_2(\w)\!=\! 0?$}}
   \put( 0,20){\line(3,-2){30}}  \put(25,7){Y}
   \put(30, 0){\line(3,2){30}}
   \put(60,20){\line(1,1){20}}   \put(63,37){N}
   \put(60,60){\line(1,-1){20}}
   \end{picture}
   }

   \put(280,300){\framebox(80,20)[c]{$\#X^2 = 6$}}

   \put(324,293){\vector(0,-1){24}}  \put(326,283){${}_6$}
   \put(330,200){${}_2$}
   \put(324,203){\line(0,-1){13}}
   \put(324,190){\vector(-1,0){15}}
   \put(285,180){\discard}

   \put(290,206){
   \begin{picture}(60,60)(0,0)
   \put( 0,20){\framebox(60,40)[c]{$f_1(\w)\!=\! 0?$}}
   \put( 0,20){\line(3,-2){30}}  \put(25,7){Y}
   \put(30, 0){\line(3,2){30}}
   \put(0,20){\line(-1,1){20}}   \put(-12,37){N}
   \put(0,60){\line(-1,-1){20}}
   \end{picture}
   }

   \put(87,245){\vector(1,0){35}}    \put(93,250){${}_3$}
   \put(270,245){\vector(-1,0){70}}  \put(260,250){${}_4$}
   \put(155,205){\vector(-1,-1){28}} \put(143,200){${}_7$}
   \put(165,205){\line(1,-1){15}}    \put(174,200){${}_5$}
   \put(180,190){\vector(1,0){15}}
   \put(195,180){\discard}

   \put(160,240){\oval(60,60)}
   \put(130,235){\begin{tabular}{c}
                   diagonal \\ homotopy \\ $3 \times 4$
                \end{tabular}}

   \put(  0,150){\framebox(80,20)[c]{$\#W^2_1 = 2$}}
   \put( 80,150){\framebox(80,20)[c]{$\#W^2_2 = 7$}}
   \put(160,150){\framebox(80,20)[c]{$W^2_3 = \emptyset$}}

   \put(33,144){\vector(-1,-1){24}}  \put(30,133){${}_2$}
   \put(13,53){\vector(1,-1){26}}    \put(20,50){${}_2$}
   \put(114,144){\vector(0,-1){24}}  \put(117,138){${}_7$}
   \put(114,53){\vector(0,-1){26}}   \put(117,50){${}_6$}

   \put(324,144){\vector(0,-1){24}}   \put(327,138){${}_8$}
   \put(327,50){${}_7$}
   \put(324,53){\vector(0,-1){20}}
   \put(312,10){\discard}

   \put(280,150){\framebox(80,20)[c]{$\#X^3 = 8$}}

   \put(290,56){
   \begin{picture}(60,60)(0,0)
   \put( 0,20){\framebox(60,40)[c]{$\begin{array}{c}
                                       f_1(\w)\!=\! 0? \\
                                       f_2(\w)\!=\! 0?
                                    \end{array}$}}
   \put( 0,20){\line(3,-2){30}}  \put(25,7){Y}
   \put(30, 0){\line(3,2){30}}
   \put(0,20){\line(-1,1){20}}   \put(-12,37){N}
   \put(0,60){\line(-1,-1){20}}
\end{picture}
}

\put(-20,56){
\begin{picture}(60,60)(0,0)
\put( 0,20){\framebox(60,40)[c]{$f_3(\w)\!=\! 0?$}}
\put( 0,20){\line(3,-2){30}}  \put(25,7){Y}
\put(30, 0){\line(3,2){30}}
\put(60,20){\line(1,1){20}}   \put(63,37){N}
\put(60,60){\line(1,-1){20}}
\end{picture}
}

\put(80,56){
\begin{picture}(60,60)(0,0)
\put( 0,20){\framebox(60,40)[c]{$f_3(\w)\!=\! 0?$}}
\put( 0,20){\line(3,-2){30}}  \put(25,7){Y}
\put(30, 0){\line(3,2){30}}
\put(60,20){\line(1,1){20}}   \put(63,37){N}
\put(60,60){\line(1,-1){20}}
\end{picture}
}

\put(167,95){\vector(1,0){18}}   \put(171,100){${}_1$}
\put(272,95){\vector(-1,0){18}}  \put(265,100){${}_1$}
\put(220,55){\vector(-1,-1){26}} \put(222,51){${}_1$}

\put(220,90){\oval(60,60)}
\put(190,85){\begin{tabular}{c}
                diagonal \\ homotopy \\ $1 \times 1$
             \end{tabular}}

\put(  0,0){\framebox(80,20)[c]{$\#W^3_1 = 2$}}
\put( 80,0){\framebox(80,20)[c]{$\#W^3_2 = 6$}}
\put(160,0){\framebox(80,20)[c]{$\#W^3_3 = 1$}}


\end{picture}
\caption{Flowchart for the illustrative example.}
\label{Fig:Flowillusex}
\end{center}
\end{figure}
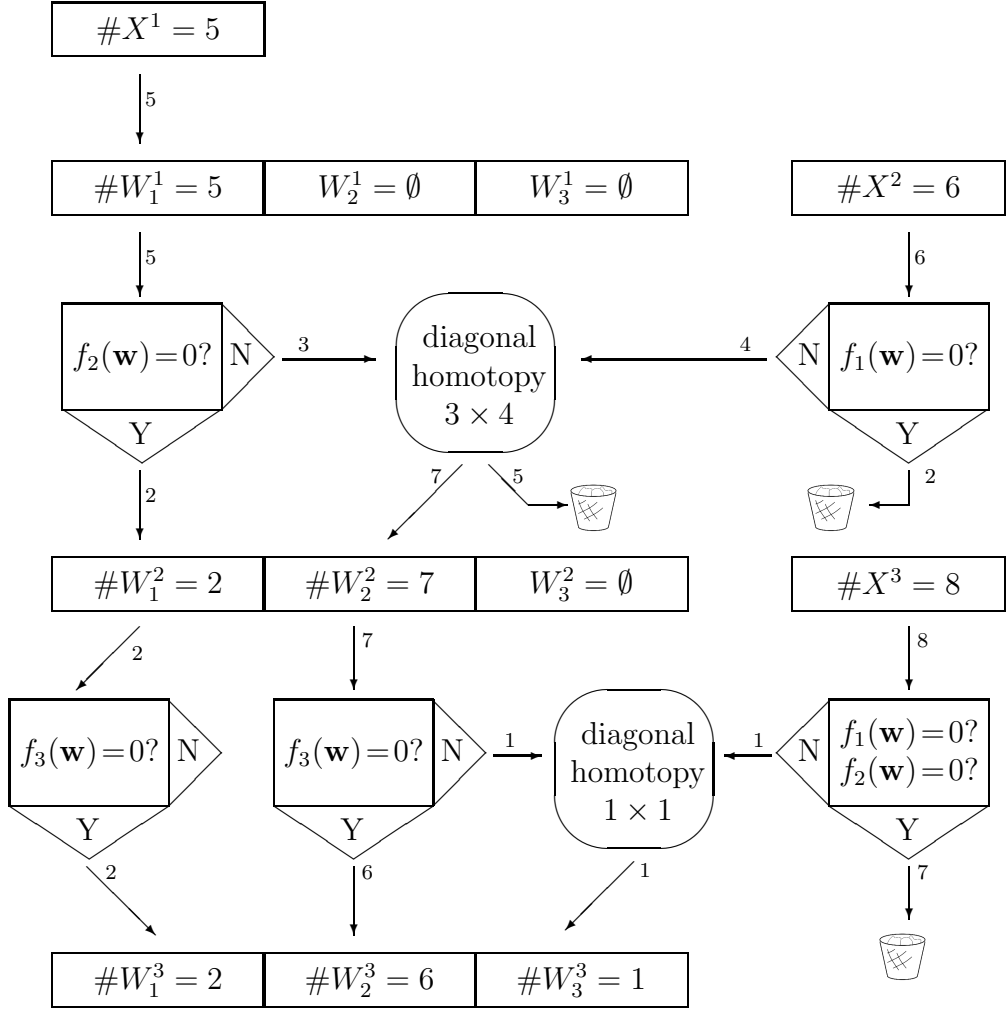

The illustrative example (see Eq.~\ref{Eq:Illusex} for the system)
illustrates the gains made by our new solver. While our previous
sequence of homotopies needed 197 paths to find all candidate
witness points, the new approach shown in
Figure~\ref{Fig:Flowillusex} tracks just 13 paths.  Many of the paths
take shortcuts around the diagonal homotopies, and five paths that
diverge to infinity in the first diagonal homotopy need no further
consideration.  It happens that none of the witness points generated
by the diagonal homotopies is singular, so there is no need for
membership testing.

On a 2.4Ghz Linux workstation, our previous approach~\cite{SV}
requires a total of 43.3 cpu seconds (39.9 cpu seconds for solving
the top dimensional embedding and 3.4 cpu seconds to run the cascade
of homotopies to find all candidate witness points). Our new
approach takes slightly less than a second of cpu time. So for this
example our new solver is 40 times faster.

\subsection{Adjacent Minors of a General 2-by-9 Matrix}

In an application from algebraic statistics~\cite{DES98}
(see also~\cite{HS00} for methods dedicated for these type of ideals)
one considers all adjacent minors of a general matrix.
For instance, consider this general 2-by-9 matrix:
\begin{displaymath}
  \left[
     \begin{array}{ccccccccc}
        x_{11} & x_{12} & x_{13} & x_{14} &
        x_{15} & x_{16} & x_{17} & x_{18} & x_{19} \\
        x_{21} & x_{22} & x_{23} & x_{24} &
        x_{25} & x_{26} & x_{27} & x_{28} & x_{29} \\
   \end{array}
 \right]
\end{displaymath}
Two minors are adjacent if they share one neighboring column.
Taking all adjacent minors from this general 2-by-9 matrix gives
8 quadrics in 18 unknowns.  This defines a 10-dimensional surface,
of degree 256.

We include this example to illustrate that the flow of timings is
typical as in Table~\ref{tabminors}. Although we execute many
homotopies, most of the work occurs in the last stage, because both
the number of paths and the number of variables increases at each
stage.  We are using the intrinsic method of \cite{SVW10} to reduce
the number of variables.  With the older extrinsic method of
\cite{SVW9}, the total cpu time increases five-fold from 104s
to~502s.

\begin{table}[hbt]
\begin{center}
\begin{tabular}{|c|rcr|c|r|} \hline
  stage & \multicolumn{3}{c|}{\#paths}
        & time/path & time \\ \hline
    1   &   4 & = &   2 $\times$ 2 & 0.03s &   0.11s \\
    2   &   8 & = &   4 $\times$ 2 & 0.05s &   0.41s \\
    3   &  16 & = &   8 $\times$ 2 & 0.10s &   1.61s \\
    4   &  32 & = &  16 $\times$ 2 & 0.12s &   3.75s \\
    5   &  64 & = &  32 $\times$ 2 & 0.19s &  12.41s \\
    6   & 128 & = &  64 $\times$ 2 & 0.27s &  34.89s \\
    7   & 256 & = & 128 $\times$ 2 & 0.41s & 104.22s \\ \hline
  \multicolumn{5}{|r}{total user cpu time\quad} & 157.56s \\ \hline
\end{tabular}
\caption{Timings on Apple PowerBook G4 1GHz for the $2\times9$
adjacent minors, a system of 8 quadrics in 18 unknowns.}
\label{tabminors}
\end{center}
\end{table}

\subsection{A General 6-by-6 Eigenvalue Problem}

Consider $f(\x,\lambda) = \lambda \x - A \x = \zero$, where $A \in
\bC^{6 \times 6}$, $A$ is a random matrix. These 6 equations in 7
unknowns define a curve of degree~7, far less than what may be
expected from the application of B\'ezout's theorem: $2^6=64$.
Regarded as a polynomial system on $\bC^7$, the solution set
consists of seven lines, six of which are eigenvalue-eigenvector
pairs while the seventh is the trivial line $\x = \zero$.

Clearly, as a matter of practical computation, one would employ an
off-the-shelf eigenvalue routine to solve this problem efficiently.
Even with continuation, we could cast the problem on
$\pn1\times\pn6$ and solve it with a seven-path two-homogeneous
formulation.  However, for the sake of illustration, let us consider
how the equation-by-equation approach performs, keeping in mind that
the only information we use about the structure of the system is the
degree of each equation.  That is, we treat it just like any other
system of 6 quadratics in 7 variables and let the
equation-by-equation procedure numerically discover its special
structure.

In a direct approach of solving the system in one total-degree
homotopy, adding one generic linear equation to slice out an
isolated point on each solution line, we would have 64 paths of
which 57 diverge. This does not even consider the work that would be
needed if we wanted to rigorously check for higher dimensional
solution sets.

Table~\ref{tabeigen} shows the evolution of the number of solution
paths tracked in each stage of the equation-by-equation approach.
The size of each initial witness set is $\#(X^i)=2$, so each new
stage tracks two paths for every convergent path in the previous
stage.  If the quadratics were general, this would build up
exponentially to 64 paths to track in the final stage, but the
special structure of the eigenvalue equations causes there to be
only $i+2$ solutions at the end of stage~$i$. Accordingly, there are
only 12 paths to track in the final, most expensive stage, and only
40 paths tracked altogether. The seven convergent paths in the final
stage give one witness point on each of the seven solution lines.

\begin{table}[hbt]
\begin{center}
  \begin{tabular}{|c|rrrrr|c|} \hline
      stage in solver  & ~1 & ~2 & ~3 & ~4 & ~5 & total \\ \hline
    \#paths tracked    & ~4 & ~6 & ~8 & 10 & 12 &  40 \\ \hline
    \#divergent paths  & ~1 & ~2 & ~3 & ~4 & ~5 &  15 \\
    \#convergent paths & ~3 & ~4 & ~5 & ~6 & ~\raise1pt\hbox to0pt{\hskip-3pt$\bigcirc$\hss}{7} &  25 \\ \hline
  \end{tabular}
\caption{Number of convergent and divergent paths on a general
         6-by-6 eigenvalue problem.}
\label{tabeigen}
\end{center}
\end{table}

\section{Conclusions}

The recent invention of the diagonal homotopy allows one to compute
intersections between algebraic sets represented numerically by
witness sets.  This opens up many new possibilities for ways to
manipulate algebraic sets numerically.  In particular, one may solve
a system of polynomial equations by first solving subsets of the
equations and then intersecting the results.  We have presented a
subsystem-by-subsystem algorithm based on this idea, which when
carried to extreme gives an equation-by-equation algorithm.  The
approach can generate witness sets for all the solution components
of a system, or it can be specialized to only seek the nonsingular
solutions at the lowest dimension.  Applying this latter form to a
system of $N$ equations in $N$ variables, we come full circle in the
sense that we are using methods developed to deal with
higher dimensional solution sets as a means of finding just the
isolated solutions.

Experiments with a few simple systems indicates that the method can
be very effective.  Using only the total degrees of the equations,
the method numerically discovers some of their inherent structure in
the early stages of the computation.  These early stages are
relatively cheap and they can sometimes eliminate much of the
computation that would otherwise be incurred in the final stages.

In future work, we plan to exercise the approach on more challenging
problems, especially ones where the equations have
interrelationships that are not easily revealed just by examining
the monomials that appear.  Multihomogenous homotopies and
polyhedral homotopies are only able to take advantage of that sort
of structure, while the equation-by-equation approach can reveal
structure encoded in the coefficients of the polynomials. One avenue
of further research could be to seek a formulation that uses
multihomogeneous homotopies or polyhedral homotopies in an
equation-by-equation style to get the best of both worlds.

\end{document}